\RequirePackage{ifpdf}
\ifpdf 
\documentclass[pdftex]{sigma}
\else
\documentclass{sigma}
\fi

\newcommand{\pa}{{\partial}}

\newcommand{\g}{{\bf g}}
\newcommand{\pxi}{{\pa \over \pa x^i}}

\begin{document}

\allowdisplaybreaks

\renewcommand{\PaperNumber}{008}

\FirstPageHeading

\ShortArticleName{On Special Berwald Metrics}

\ArticleName{On Special Berwald Metrics}

\Author{Akbar TAYEBI~$^\dag$ and Esmaeil PEYGHAN~$^\ddag$}

\AuthorNameForHeading{A.~Tayebi and E.~Peyghan}

\Address{$^\dag$~Department of Mathematics, Faculty  of Science, Qom University, Qom, Iran}
\EmailD{\href{mailto:akbar.tayebi@gmail.com}{akbar.tayebi@gmail.com}}

\Address{$^\ddag$~Department of Mathematics, Faculty  of Science, Arak University, Arak, Iran}
\EmailD{\href{mailto:epeyghan@gmail.com}{epeyghan@gmail.com}, \href{e-peyghan@araku.ac.ir}{e-peyghan@araku.ac.ir}}

\ArticleDates{Received November 01, 2009, in f\/inal form January 17, 2010;  Published online January 20, 2010}

\Abstract{In this paper, we study  a class of Finsler metrics which contains the class of Berwald metrics as a special case. We prove that every Finsler metric in this class is a~generalized Douglas--Weyl metric.  Then we study isotropic  f\/lag curvature Finsler metrics in this class.  Finally we show that on this class of Finsler metrics, the notion of Landsberg and weakly Landsberg curvature are equivalent.}


\Classification{53C60; 53C25}

\section{Introduction}

For a Finsler metric $F=F(x,y)$, its geodesics curves are
characterized by the system of dif\/fe\-ren\-tial equations $ \ddot
c^i+2G^i(\dot c)=0$, where the local functions $G^i=G^i(x, y)$
are called the spray coef\/f\/icients. A Finsler metric $F$ is
called a Berwald metric if $G^i = {1\over 2}
\Gamma^i_{jk}(x)y^jy^k$ is quadratic in $y\in T_xM$  for any
$x\in M$. It is  proved that on  a Berwald space,  the parallel
translation along any geodesic preserves the Minkowski
functionals \cite{Ich}. Thus Berwald spaces can be viewed as
Finsler spaces modeled on a single Minkowski space.

Recently by using the structure of Funk metric, Chen--Shen
introduce the notion of isotropic Berwald metrics
\cite{ChSh,TR}. This motivates us to study special forms
of Berwald metrics.

Let $(M, F)$ be a two-dimensional Finsler manifold. We refer to
the Berwald's frame $(\ell^i, m^i)$ where $\ell^i=y^i/F(y)$,
$m^i$ is the unit vector with $\ell_i m^i=0$,
$\ell_i=g_{ij}\ell^i$ and $g_{ij}$ is the fundamental tensor of
Finsler metric $F$.  Then the Berwald curvature is given by
\begin{gather}\label{TP1}
B^i_{\ jkl}=F^{-1}\big({-}2I_{,1}\ell^i+I_2m^i\big) m_j m_k m_l,
\end{gather}
where $I$ is 0-homogeneous function called the main scalar of
Finsler metric and $I_2=I_{,2}+I_{,1|2}$ (see~\cite[page~689]{An}). By (\ref{TP1}), we have
\begin{gather*}
B^i_{\
jkl}=-\frac{2I_{,1}}{3F^2}\big(m_jh_{kl}+m_kh_{jl}+m_lh_{jk}\big)y^i+\frac{I_2}{3F}\big(h^i_j
h_{kl}+h^i_k h_{jl}+h^i_l h_{jk}\big),
\end{gather*}
where $h_{ij}:=m_i m_j$ is called the angular metric. Using the
special form of Berwald curvature for Finsler surfaces, we def\/ine
a new class of Finsler metrics on $n$-dimensional Finsler manifolds
which their Berwald curvature satisfy in following
\begin{gather}\label{TP3}
B^i_{\
jkl}=(\mu_jh_{kl}+\mu_kh_{jl}+\mu_lh_{jk})y^i+\lambda\big(h^i_j
h_{kl}+h^i_k h_{jl}+h^i_l h_{jk}\big),
\end{gather}
where $\mu_i=\mu_i(x,y)$ and $\lambda=\lambda(x,y)$ are
homogeneous functions of degrees~$-2$ and~$-1$ with~respect to~$y$,
respectively. By def\/inition of Berwald curvature, the function $\mu_i$ satisf\/ies \mbox{$\mu_i y^i{=}0$}~\cite{PTP}.

The Douglas tensor is another non-Riemanian curvature def\/ined as
follows
\begin{gather}\label{TP11}
D^i_{\ jkl}:=\left(G^i-\frac{1}{n+1}\frac{\partial G^m}{\partial
y^m}y^i\right)_{y^j y^k y^l}.
\end{gather}
Douglas curvature is a non-Riemannian projective invariant
constructed from the Berwald curvature. The notion of Douglas
curvature was proposed by B\'acs\'o and Matsumoto
as a gene\-ra\-lization of Berwald curvature \cite{BM}. We show that
a Finsler metric satisf\/ies (\ref{TP3}) with vanishing Douglas
tensor is a Randers metric (see Proposition \ref{Prop1}). A
Finsler metric is called a genera\-li\-zed Douglas--Weyl (GDW) metric
if the Douglas tensor satisfy in $h^i_{\alpha} D^\alpha_{\
jkl|m}y^m=0$~\cite{NST1}. In \cite{BP}, B\'{a}cs\'{o}--Papp show
that this class of Finsler metrics is closed under projective
transformation. We prove that a Finsler metric satisf\/ies
(\ref{TP3}) is a GDW-metric.

\begin{theorem}\label{THM1}
Every  Finsler metric satisfying  \eqref{TP3} is a GDW-metric.
\end{theorem}
Theorem \ref{THM1}, shows that every two-dimensional Finsler
metric is a generalized Douglas--Weyl metric.

For a Finsler manifold $(M, F)$, the f\/lag curvature is  a
function ${\bf K}(P, y)$ of tangent planes $P\subset T_xM$ and
directions $y\in P$. $F$  is said to be of  isotropic f\/lag
curvature if  ${\bf K}(P, y)={\bf K}(x)$  and constant f\/lag
curvature if ${\bf K}(P, y)={\rm const}$.

\begin{theorem}\label{THM2}
Let $F$ be a Finsler metric of non-zero isotropic flag curvature
${\bf K}={\bf K}(x)$ on a~manifold $M$.  Suppose that $F$
satisfies  \eqref{TP3}. Then $F$ is a Riemannian metric if and
only if $\mu_i$ is constant along geodesics.
\end{theorem}

Beside the Berwald curvature, there are several important
Finslerian curvature. Let $(M, F)$ be a Finsler manifold. The
second derivatives of ${1\over 2} F_x^2$ at $y\in T_xM_{0}$ is an
inner product~${\bf g}_y$ on~$T_xM$.  The third order derivatives
of ${1\over 2} F_x^2$ at  $y\in T_xM_0$ is a symmetric trilinear
forms ${\bf C}_y$ on $T_xM$. We call ${\bf g}_y$ and ${\bf C}_y$
the  fundamental form and  the  Cartan torsion, respectively.
The rate of change of the Cartan torsion along geodesics is the
Landsberg curvature  ${\bf L}_y$ on $T_xM$ for any $y\in T_xM_0$.
Set ${\bf J}_y:= \sum\limits_{i=1}^n {\bf L}_y(e_i, e_i, \cdot )$, where
$\{e_i\}$ is an orthonormal basis for $(T_xM, g_y)$.  ${\bf J}_y$~is called the  mean Landsberg curvature. $F$~is said to be
Landsbergian if ${\bf L}=0$,  and  weakly Landsbergian if ${\bf
J}=0$ \cite{ShLec,ShDiff}.

In this paper, we prove that on Finsler manifolds satisf\/ies
(\ref{TP3}), the notions of  Landsberg  and  weakly Landsberg
metric are equivalent.

\begin{theorem}\label{THM3}
Let $(M, F)$ be a Finsler manifold  satisfying \eqref{TP3}.  Then
${\bf L}=0$ if and only if~${\bf J}=0$.
\end{theorem}

There are many connections in Finsler geometry  \cite{TAE}. In
this paper, we use the Berwald connection  and the $h$- and $v$-covariant derivatives of a Finsler tensor f\/ield are denoted by ``$|$'' and ``,'' respectively.

\section{Preliminaries}

Let $M$ be a n-dimensional $ C^\infty$ manifold. Denote by $T_x M
$ the tangent space at $x \in M$,  by $TM=\cup _{x \in M} T_x M $
the tangent bundle of $M$, and by $TM_{0} = TM \setminus \{ 0 \}$
the slit tangent bundle on~$M$. A  Finsler metric on $M$ is a
function $ F:TM \rightarrow [0,\infty)$ which has the following
properties: $(i)$~$F$ is $C^\infty$ on $TM_{0}$; $(ii)$~$F$ is
positively 1-homogeneous on the f\/ibers of tangent bundle $TM$,
and $(iii)$~for each $y\in T_xM$, the following quadratic form
$g_y$ on $T_xM$  is positive def\/inite,
\[
g_{y}(u,v):={1 \over 2} \frac{d^2}{ds dt}\left[  F^2
(y+su+tv)\right]\big|_{s,t=0}, \qquad u,v\in T_xM.
\]
Let  $x\in M$ and $F_x:=F|_{T_xM}$.  To measure the non-Euclidean
feature of $F_x$, def\/ine ${\bf C}_y$: $T_xM\times T_xM\times
T_xM\rightarrow \mathbb{R}$ by
\[
{\bf C}_{y}(u,v,w):={1 \over 2} \frac{d}{dt}\left[g_{y+tw}(u,v)
\right]|_{t=0}, \qquad u,v,w\in T_xM.
\]
The family $\bf{C}:=\{\bf{C}_y\}_{y\in TM_0}$  is called the
Cartan torsion. It is well known that ${\bf{C}}=0$ if and only if
$F$ is Riemannian~\cite{ShDiff}. For $y\in T_x M_0$, def\/ine  mean
Cartan torsion ${\bf I}_y$ by ${\bf I}_y(u):=I_i(y)u^i$, where
$I_i:=g^{jk}C_{ijk}$, $g^{jk}$ is the inverse of $g_{jk}$ and
$u=u^i\frac{\partial}{\partial x^i}|_x$. By Deicke's  theorem,
$F$ is Riemannian  if and only if ${\bf I}_y=0$~\cite{ShLec}.

Let $\alpha=\sqrt{a_{ij}(x)y^iy^j}$ be a Riemannian metric, and
$\beta=b_i(x)y^i$ be a 1-form on $M$ with
$b=\sqrt{a^{ij}b_ib_j}<1$. The Finsler metric $F=\alpha+\beta$ is
called a Randers metric.

Let $(M, F)$ be a Finsler manifold. Then for a non-zero vector
$y \in T_xM_0$, def\/ine the  Matsumoto torsion ${\bf
M}_y:T_xM\otimes T_xM \otimes T_xM \rightarrow \mathbb{R}$ by
${\bf M}_y(u,v,w):=M_{ijk}(y)u^iv^jw^k$ where
\[
M_{ijk}:=C_{ijk} - \tfrac{1}{n+1}  \{ I_i h_{jk} + I_j h_{ik} + I_k
h_{ij} \},\label{Matsumoto}
\]
$h_{ij}:=FF_{y^iy^j}=g_{ij}-\frac{1}{F^2}g_{ip}y^pg_{jq}y^q$ is
the angular metric and $I_i:=g^{jk}C_{ijk}$ is the mean Cartan
torsion. By def\/inition, we have  $h_{ij}y^i=0$,
$h^i_j=\delta^i_j-F^{-2}y^iy_j$, $y_j=g_{ij}y^i$, $h^i_j
h_{ik}=h_{jk}$ and $h^i_i=n-1$.  A Finsler metric $F$ is said to
be $C$-reducible if ${\bf M}_y=0$. This quantity is introduced by
Matsumoto~\cite{Mat}. Matsumoto proves that every Randers metric
satisf\/ies that ${\bf M}_y=0$. Later on, Matsumoto--H\={o}j\={o}
proves that the converse is true too.

\begin{lemma}[\cite{MH}]\label{MaHo}
A Finsler metric $F$ on a manifold of dimension $n\geq 3$
is a Randers metric if and only if  ${\bf M}_y =0$, $\forall\,
y\in TM_0$.
\end{lemma}

Let us consider the pull-back tangent bundle $\pi^*TM$ over
$TM_0$ def\/ined by
\[ \pi^*TM=\left\{(u, v)\in TM_0 \times TM_0 |\ \pi(u)=\pi(v)\right\}.\]
Let $\nabla$  be the Berwald connection.
Let $\{e_i\}^n _{i=1}$ be a  local orthonormal (with respect to
$g$) frame f\/ield for the pulled-back bundle $\pi ^* TM$ such that
$e_n=\ell$, where $\ell$ is the canonical section   of $\pi^*TM$
def\/ined by $\ell_y=y/F(y)$. Let $\{\omega^i\}^n_{i=1}$ be its
dual co-frame f\/ield.  Put   $\nabla e_i = \omega^{j}_{\ i}
\otimes e_j$,  where $\{ \omega^{j}_{\ i}\}$  is called the
connection forms of $\nabla$ with respect to   $\{e_{i}\}$. Put
$\omega^{n+i}:=\omega^{i}_{\ n} +d(\log F)\delta^i_n$.  It is easy
to show that $\{\omega^i, \omega^{n+i} \}^n_{i=1}$ is a  local
basis for $T^*( TM_0)$. Since $\{\Omega^{j}_{\ i}\}$ are 2-forms
on $TM_0$,  they can be expanded as
\[
\Omega^j_{\ i}=\tfrac{1}{2}R^{j} _{\ ikl} \omega ^k \wedge
\omega^l +B^{j} _{\ ikl} \omega ^k \wedge \omega^{n+l}.
\]
Let $\{\bar e_i, \dot e_i\}^n _{i=1}$ be the local basis for
$T(TM_0)$, which is dual to $\{\omega ^i, \omega^{n+i} \}^n
_{i=1}$. The objects $R$ and~$B$ are called, respectively,  the
$hh$- and $hv$-curvature  tensors of the Berwald connection  with the
components $R(\bar e_k,\bar e_l)e_i =R^{j}_{\ ikl}e_j$ and $
P(\bar e_k,\dot e_l)e_i=P^{j}_{\ ikl} e_j$~\cite{TAE}.  With
the Berwald connection, we def\/ine covariant derivatives of
quantities on $TM_0$  in the usual way. For example, for a scalar
function $f$, we def\/ine $f_{|i}$ and $ f_{\cdot i}$ by
\[
df = f_{|i} \omega^i + f_{, i} \omega^{n+i},
\]
where ``$|$'' and ``,''  denote the $h$- and $v$-covariant
derivatives, respectively.

The horizontal covariant derivatives of ${\bf C}$ along geodesics
give rise to  the  Landsberg curvature  ${\bf L}_y:T_xM\times
T_xM\times T_xM\rightarrow \mathbb{R}$  def\/ined by
\[
{\bf L}_y(u,v,w):=L_{ijk}(y)u^iv^jw^k,
\]
where $L_{ijk}:=C_{ijk|s}y^s$, $u=u^i{{\partial } \over {\partial
x^i}}|_x$,  $v=v^i{{\partial }\over {\partial x^i}}|_x$ and
$w=w^i{{\partial }\over {\partial x^i}}|_x$. The family ${\bf
L}:=\{{\bf L}_y\}_{y\in TM_{0}}$  is called the \textit{Landsberg
curvature}. A Finsler metric is called a \textit{Landsberg
metric}  if {\bf{L}=0}. The horizontal covariant derivatives of
${\bf I}$ along geodesics give rise to  the mean Landsberg
curvature ${\bf J}_y(u): = J_i (y)u^i$, where $J_i: =
g^{jk}L_{ijk}$. A Finsler metric is said to be {\it weakly
Landsbergian} if ${\bf J}=0$.

Given a Finsler manifold $(M,F)$, then a global vector f\/ield $G$
is induced by $F$ on $TM_0$, which in a standard coordinate
$(x^i,y^i)$ for $TM_0$ is given by
\[
G=y^i {{\partial} \over {\partial x^i}}-2G^i(x,y){{\partial} \over
{\partial y^i}},
\]
where $G^i(y)$ are local functions on $TM$ given by
\[
G^i(y):=\frac{1}{4}g^{il}(y)\left\{\frac{\partial^2[F^2]}{\partial x^k
\partial y^l}(y)y^k-\frac{\partial[F^2]}{\partial x^l}(y)\right\},\qquad
y\in T_xM.
\]
$G$ is called the  spray associated  to $(M,F)$.  In local
coordinates, a curve $c(t)$ is a geodesic if and only if its
coordinates $(c^i(t))$ satisfy $ \ddot c^i+2G^i(\dot c)=0$.

For a tangent vector $y \in T_xM_0$, def\/ine ${\bf
B}_y:T_xM\otimes T_xM \otimes T_xM\rightarrow T_xM$ and ${\bf
E}_y:T_xM \otimes T_xM\rightarrow \mathbb{R}$ by ${\bf B}_y(u, v,
w):=B^i_{\ jkl}(y)u^jv^kw^l{{\partial } \over {\partial x^i}}|_x$
and ${\bf E}_y(u,v):=E_{jk}(y)u^jv^k$ where
\[
B^i_{\ jkl}(y):={{\partial^3 G^i} \over {\partial y^j \partial
y^k \partial y^l}}(y),\qquad E_{jk}(y):=\tfrac{1}{2} B^m_{\
jkm}(y).
\]
$\bf B$ and $\bf E$ are called the Berwald curvature
and mean Berwald curvature, respectively.  Then~$F$ is called a
Berwald metric and weakly Berwald metric if $\bf{B}=0$ and
$\bf{E}=0$, respectively~\cite{ShDiff}. By def\/inition of Berwald
and mean Berwald curvatures, we have
\begin{gather*}
y^j B^i_{\ jkl}=y^k B^i_{\ jkl}=y^l B^i_{\ jkl}=0, \qquad y^j
E_{jk}=y^k E_{jk}=0.
\end{gather*}
The Riemann curvature ${\bf R}_y= R^i_{\ k}  dx^k \otimes \pxi|_x
: T_xM \to T_xM$ is a family of linear maps on tangent spaces,
def\/ined by
\begin{gather*}
R^i_{\ k} = 2 {\pa G^i\over \pa x^k}-y^j{\pa^2 G^i\over \pa
x^j\pa y^k} +2G^j {\pa^2 G^i \over \pa y^j \pa y^k} - {\pa G^i
\over \pa y^j} {\pa G^j \over \pa y^k}.
\end{gather*}
The f\/lag curvature in Finsler geometry is a natural extension of
the sectional curvature in Riemannian geometry was  f\/irst
introduced by L.~Berwald \cite{Be}. For a f\/lag $P={\rm span}\{y,
u\} \subset T_xM$ with f\/lagpole $y$, the  f\/lag curvature ${\bf
K}={\bf K}(P, y)$ is def\/ined by
\begin{gather*}
{\bf K}(P, y):= {\g_y (u, {\bf R}_y(u)) \over \g_y(y, y) \g_y(u,u)
-\g_y(y, u)^2 }.
\end{gather*}
When $F$ is Riemannian, ${\bf K}={\bf  K}(P)$ is independent of
$y\in P$,  and is the sectional curvature of~$P$. We say that a
Finsler metric $F$ is   of scalar curvature if for any $y\in
T_xM$, the f\/lag curvature ${\bf K}= {\bf K}(x, y)$ is a scalar
function on the slit tangent bundle $TM_0$. If ${\bf
K}={\rm const}$, then $F$ is said to be of  constant f\/lag curvature.
A Finsler metric $F$ is called {\it isotropic flag curvature}, if
${\bf K}= {\bf K}(x)$.

In \cite{AZ}, Akbar-Zadeh considered a non-Riemannian quantity
$\bf{H}$ which is obtained from the mean Berwald curvature by the
covariant horizontal dif\/ferentiation along geodesics. This is a~positively homogeneous scalar function of degree zero on the slit
tangent bundle. The quantity ${\bf H}_y = H_{ij}dx^i\otimes dx^j$
is def\/ined as the covariant derivative of ${\bf E}$ along
geodesics~\cite{NST2}. More precisely
\[
H_{ij}:= E_{ij|m} y^m.
\]
In local coordinates, we have
\begin{gather*}
2H_{ij}= y^m\frac{\partial^4 G^k}{\partial y^i\partial y^j \partial y^k \partial x^m}-2G^m \frac{\partial^4 G^k}{\partial y^i \partial y^j\partial y^k\partial y^m}
- \frac{\partial G^m}{\partial
y^i}\frac{\partial^3 G^k}{\partial y^j\partial y^k\partial
y^m}-\frac{\partial G^m}{\partial y^j} \frac{\partial^4
G^k}{\partial y^i\partial y^k\partial y^m}.
\end{gather*}
Akbar-Zadeh proved the following:
\begin{theorem}[\cite{AZ}]\label{Akbar}   Let $F$ be a Finsler metric of scalar curvature on an n-dimensional manifold~$M$ $(n\geq 3)$. Then  the flag curvature ${\bf K}={\rm const}$  if and only if~${\bf H}=0$.
\end{theorem}

\section{Proof of Theorem \ref{THM1}}

\begin{lemma}\label{Lem1}
Let $(M,F)$ be a Finsler manifold. Suppose that the Cartan tensor
satisfies in $C_{ijk}=B_ih_{jk}+B_jh_{ik}+B_kh_{ij}$ with
$y^iB_i=0$. Then $F$ is a $C$-reducible metric.
\end{lemma}

\begin{proof}
Suppose that the Cartan tensor of the Finsler metric $F$
satisf\/ies in
\begin{gather}\label{TP6}
C_{ijk}=B_ih_{jk}+B_jh_{ik}+B_kh_{ij}.
\end{gather}
Contracting (\ref{TP6}) with $g^{ij}$ yields
\begin{gather}\label{TP7}
I_k=B_ih^i_k+B_jh^j_k+(n-1)B_k.
\end{gather}
Using (\ref{TP7}) and $B_ih^i_k=B_jh^j_k=B_k$, we get
$I_i=(n+1)B_i$. Putting this relation in (\ref{TP6}), we conclude
that $F$  is a $C$-reducible Finsler metric.
\end{proof}

\begin{lemma}\label{Lem2}
Let $(M,F)$ be a Finsler metric. Then $F$ is a GDW-metric if and
only if
\begin{gather}\label{TP8}
D^i_{\ jkl|s}y^s=T_{jkl}y^i,
\end{gather}
for some tensor $T_{jkl}$ on manifold $M$.
\end{lemma}
\begin{proof}
Let $F$ be is a GDW-metric
\begin{gather*}
h^i_mD^m_{\ jkl|s}y^s=0.
\end{gather*}
This yields
\begin{gather*}
D^i_{\ jkl|s}y^s=\big(F^{-2} y_mD^m_{\ jkl|s}\big)y^i.
\end{gather*}
Therefore $T_{jkl}:=F^{-2} y_mD^m_{\ jkl|s}$. The proof of
converse is trivial.
\end{proof}

Equation (\ref{TP8})  is equivalent to the condition that, for
any  parallel vector f\/ields $U=U(t)$, $V=V(t)$ and $W=W(t)$ along
a geodesic $c(t)$, there is a function $T=T(t)$ such that
\[
\frac{d}{dt}[D_{\dot{c}}(U,V,W)]=T\dot{c}.
\]
The geometric meaning of the above identity is that the rate of
change of the Douglas curvature along a geodesic is tangent to
the geodesic.

\begin{proposition}\label{Prop1}
Let $(M,F)$ be a Finsler manifold satisfies \eqref{TP3} with
dimension $n\geq 3$. Suppose that the Douglas tensor of $F$
vanishes. Then $F$ is a Randers metric.
\end{proposition}
\begin{proof}
Since $F$ satisf\/ies  (\ref{TP3}), then  by considering $\mu_i
y^i=0$ we get
\begin{gather}\label{TP12}
2E_{jk}=(n+1)\lambda h_{ij}.
\end{gather}
On the other hand, we have
\begin{gather*}
h_{ij,k}=2C_{ijk}-F^{-2}(y_jh_{ik}+y_ih_{jk}),
\end{gather*}
which implies that
\begin{gather}\label{TP14}
2E_{jk,l}=(n+1)\lambda_{,l}
h_{jk}+(n+1)\lambda\big\{2C_{jkl}-F^{-2}(y_kh_{jl}+y_jh_{kl})\big\}.
\end{gather}
Putting (\ref{TP3}), (\ref{TP12}) and (\ref{TP14}) in
(\ref{TP11}) yields
\begin{gather}\label{TP15}
D^i_{\ jkl}=\{\mu_jh_{kl}+\mu_kh_{jl}+\mu_lh_{jk}-2\lambda
C_{jkl}\}y^i-\big(\lambda y_lF^{-2}+\lambda_{,l}\big)h_{jk} y^i.
\end{gather}
For the Douglas curvature, we have $D^i_{\ jkl}=D^i_{\ jlk}$.
Then by (\ref{TP15}), we conclude that
\begin{gather}\label{TP16}
\lambda y_lF^{-2}+\lambda_{,l}=0.
\end{gather}
From (\ref{TP15}) and (\ref{TP16}) we deduce
\begin{gather}\label{TP17}
D^i_{\ jkl}=\{\mu_jh_{kl}+\mu_kh_{jl}+\mu_lh_{jk}-2\lambda
C_{jkl}\}y^i.
\end{gather}
Since $F$ is a Douglas metric, then
\begin{gather*}
C_{jkl}=\tfrac{1}{2\lambda}\{\mu_jh_{kl}+\mu_kh_{jl}+\mu_lh_{jk}\}.
\end{gather*}
By Lemmas \ref{Lem1} and~\ref{MaHo}, it follows that $F$ is
a Randers metric.
\end{proof}

\begin{proof}[Proof of Theorem \ref{THM1}]  To prove the
Theorem \ref{THM1}, we start with the equation (\ref{TP17}):
\begin{gather}\label{TP19}
D^i_{\ jkl}=\{\mu_jh_{kl}+\mu_kh_{jl}+\mu_lh_{jk}-2\lambda
C_{jkl}\}y^i.
\end{gather}
Taking a horizontal derivation of (\ref{TP19}) implies that
\begin{gather*}
D^i_{\
jkl|s}y^s=\{\mu'_jh_{kl}+\mu'_kh_{jl}+\mu'_lh_{jk}-2\lambda'
C_{jkl}-2\lambda L_{jkl}\}y^i.
\end{gather*}
where $\lambda'=\lambda_{|m} y^m$ and $\mu'_i=\mu_{i|m} y^m$. By
Lemma \ref{Lem2}, $F$ is a GDW-metric with
\begin{gather*}
T_{jkl}=\mu'_jh_{kl}+\mu'_kh_{jl}+\mu'_lh_{jk}-2\lambda'
C_{jkl}-2\lambda L_{jkl}.
\end{gather*}
This completes the proof.
\end{proof}

  The  Funk metric on a  strongly convex domain
$\mathbb{B}^n\subset\mathbb{R}^n$ is a non-negative function on
$T\Omega=\Omega\times \mathbb{R}^n$, which in the special case
$\Omega=\mathbb{B}^n$ (the unit ball in the Euclidean space
$\mathbb{R}^n$) is def\/ined by the following explicit formula:
\[
F(y):=\frac{\sqrt{|y|^2-(|x|^2|y|^2-\langle x,y\rangle ^2)}}{1-|x|^2}+\frac{\langle x,y\rangle }{1-|x|^2},\qquad
  y\in T_x\mathbb{B}^n=\mathbb{R}^n,
\]
where $|\cdot|$ and $\langle \cdot ,\cdot \rangle$ denote the Euclidean norm and inner product
in $\mathbb{R}^n$, respectively~\cite{ShDiff}. The Funk metric on
$\mathbb{B}^n$ is a Randers metric.  The Berwald curvature of
Funk metric is given by
\[
B^i_{\ jkl}=\tfrac{1}{2F}\big\{h^i_j h_{kl}+h^i_k h_{jl}+h^i_l
h_{jk}+2C_{jkl} y^i\big\}.
\]
Thus the Funk metric is a GDW-metric which does not satisfy
(\ref{TP3}). Then by Theorem \ref{THM1}, we conclude the
following.

\begin{corollary}
The class of Finsler metrics satisfying \eqref{TP3} is a proper
subset of the class of generalized Douglas--Weyl metrics.
\end{corollary}

\section{Proof of Theorem \ref{THM2}}

To prove Theorem \ref{THM2}, we need the following.
\begin{lemma}[\cite{Ich,NST2}]\label{Lem3}
For the Berwald connection, the following Bianchi
identities hold:
\begin{gather}
R^i_{\ jkl|m}+ R^i_{\ jlm|k}+R^i_{\ jmk|l}=0,\nonumber\\ 
B^i_{\ jml|k}- B^i_{\ jkm|l}=R^i_{\ jkl,m},\label{TP23}\\
B^i_{\ jkl,m}=B^i_{\ jkm,l}.\nonumber
\end{gather}
\end{lemma}

\begin{proof}[Proof of Theorem \ref{THM2}] We have:
\begin{gather}\label{TP25}
R^i_{\ jkl} = \frac{1}{3}\left\{\frac{\partial^2
R^i_{\ k}}{\partial y^j \partial y^l}-\frac{\partial^2 R^i_{\
l}}{\partial y^j \partial y^k}\right\}.
\end{gather}
Here, we assume that a Finsler metric $F$ is of isotropic f\/lag
curvature ${\bf K}={\bf K}(x)$. In local coordinates, $R^i_{\ k}
= \textbf{K}(x)F^2h^i_k$. Plugging this equation into
(\ref{TP25}) gives
\begin{gather}\label{TP26}
R^i_{\ jkl}={\bf K}\{g_{jl}\delta^i_k-g_{jk}\delta^i_l\}.
\end{gather}
Dif\/ferentiating (\ref{TP26}) with respect to $y^m$ gives a
formula for $R^i_{\  jkl,m}$ expressed in terms of ${\bf K}$ and
its derivatives. Contracting  (\ref{TP23}) with $y^k$, we obtain
\begin{gather}\label{TP27}
B^i_{\ jml|k}y^k= 2{\bf K}C_{jml}y^i.
\end{gather}
Multiplying  (\ref{TP27}) with $y_i$ implies that
\begin{gather}\label{TP28}
B^i_{\ jml|k}y^ky_i= 2{\bf K}F^2C_{jml}.
\end{gather}
Since $F$ satisf\/ies (\ref{TP3}), then we have
\begin{gather}\label{TP29}
B^i_{\
jkl|m}y^m=(\mu'_jh_{kl}+\mu'_kh_{jl}+\mu'_lh_{jk})y^i+\lambda'(h^i_j
h_{kl}+h^i_k h_{jl}+h^i_l h_{jk}).
\end{gather}
By contracting (\ref{TP29}) with $y_i$, we have
\begin{gather}\label{TP30}
B^i_{\ jkl|m}y^my_i=(\mu'_jh_{kl}+\mu'_kh_{jl}+\mu'_lh_{jk})F^2.
\end{gather}
By (\ref{TP28}) and (\ref{TP30}) we get
\begin{gather*}
\mu'_jh_{kl}+\mu'_kh_{jl}+\mu'_lh_{jk}=2{\bf K}C_{jkl}.
\end{gather*}
Contracting with $g^{kl}$ yields
\begin{gather*}
\mu'_j=\frac{2{\bf K}}{n+1}I_j.
\end{gather*}
Since ${\bf K}\neq 0$, then  by Deicke's theorem $F$ is a
Riemannian metric if and only if $\mu'_j=0$.
\end{proof}

\begin{theorem}
Let $F$ be a Finsler metric on an n-dimensional manifold $M$
$(n\geq 3)$ and satisfies~\eqref{TP3}. Suppose that $F$ is of
scalar flag curvature ${\bf K}$. Then ${\bf K}={\rm const}$ if and
only if~$\lambda'=0$.
\end{theorem}

\begin{proof}
Contracting $i$ and $l$ in (\ref{TP3}) yields
\begin{gather*}
2E_{jk}=(n+1)\lambda h_{jk}.
\end{gather*}
By taking a horizontal derivative  of this equation, we have
\begin{gather*}
2H_{jk}=(n+1)\lambda' h_{jk}.
\end{gather*}
Therefore $H_{jk}=0$ if and only if $\lambda'=0$. By Theorem
\ref{Akbar}, we get the proof.
\end{proof}

\section{Proof of Theorem \ref{THM3}}

In this section, we are going to prove Theorem \ref{THM3}.

\begin{proof}[Proof of Theorem \ref{THM2}] Let $F$ be a Finsler
metric satisfy in following
\begin{gather}\label{TP35}
B^i_{\
jkl}=(\mu_jh_{kl}+\mu_kh_{jl}+\mu_lh_{jk})y^i+\lambda\big(h^i_j
h_{kl}+h^i_k h_{jl}+h^i_l h_{jk}\big),
\end{gather}
where $\mu_i=\mu_i(x,y)$ and $\lambda=\lambda(x,y)$ are
homogeneous functions of degrees $-2$ and $-1$ with respect to~$y$,
respectively. Contracting (\ref{TP35}) with $y_i$ yields
\begin{gather}\label{TP39}
y_iB^i_{\ jkl}=F^2(\mu_jh_{kl}+\mu_kh_{jl}+\mu_lh_{jk})+\lambda
y_i\big(h^i_j h_{kl}+h^i_k h_{jl}+h^i_l h_{jk}\big).
\end{gather}
On the other hand, we have
\begin{gather}
  y_i B^i_{\ jkl}=-2L_{jkl},\label{TP40}\\
 y_i h^i_m=y_i\big(\delta^i_m -F^{-2} y^i y_m\big)=0.\label{TP41}
\end{gather}
See \cite[page 84]{ShDiff}. Using (\ref{TP39}), (\ref{TP40}) and
(\ref{TP41}), we get
\begin{gather}\label{TP42}
L_{jkl}=-\tfrac{1}{2}F^2\{\mu_jh_{kl}+\mu_kh_{jl}+\mu_lh_{jk}\}.
\end{gather}
By (\ref{TP42}),  it is obvious that if $\mu_i=0$ then
$L_{jkl}=0$. Conversely let $F$ be a Landsberg metric. Then we
have
\begin{gather}\label{TP43}
\mu_jh_{kl}+\mu_kh_{jl}+\mu_lh_{jk}=0.
\end{gather}
Contracting  (\ref{TP43}) with $g^{kl}$  yields $\mu_j=0$. Then
$F$ is a Landsberg metric if and only if $\mu_j=0$. Now,
contracting  (\ref{TP42}) with $g^{kl}$ yields
\begin{gather}\label{TP44}
J_j=-\tfrac{1}{2}(n+1)F^2\mu_j.
\end{gather}
By (\ref{TP44}), $J_j=0$ if and only if  $\mu_j=0$. Then ${\bf
L}=0$ if and only if ${\bf J}=0$.
\end{proof}

By using the notion of Landsberg curvature, we def\/ine the stretch curvature ${\bf \Sigma}_y:T_xM\otimes T_xM \otimes T_xM  \otimes
T_xM\rightarrow \mathbb{R}$ by ${\bf \Sigma}_y(u, v,
w,z):={\Sigma}_{ijkl}(y)u^iv^jw^kz^l$ where
\[
{\Sigma}_{ijkl}:=2(L_{ijk|l}-L_{ijl|k}).
\]
In \cite{Be}, L. Berwald has introduce the stretch curvature tensor ${\bf \Sigma}$
and showed that this tensor vanishes if and only if the length of
a vector remains unchanged under the parallel displacement  along
an inf\/initesimal parallelogram.

\begin{theorem}
Let $(M, F)$ be a Finsler manifold on which \eqref{TP3} holds.
Suppose that $F$ is a stretch metric. Then  $\mu_j$ is constant
along any Finslerian geodesics.
\end{theorem}

\begin{proof}
Taking a horizontal derivation of  (\ref{TP42}) yields
\begin{gather*}
L_{ijk|l}=-\tfrac{1}{2}F^2\{\mu_{i|l}h_{jk}+\mu_{j|l}h_{ki}+\mu_{k|l}h_{ij}\}.
\end{gather*}
Suppose that ${\bf \Sigma}=0$. Then by $L_{ijk|l}=L_{ijl|k}$, we
get
\begin{gather}\label{TP46}
\mu_{i|l}h_{jk}+\mu_{j|l}h_{ki}+\mu_{k|l}h_{ij}=\mu_{i|k}h_{jl}+\mu_{j|k}h_{li}+\mu_{l|k}h_{ij}.
\end{gather}
Multiplying (\ref{TP46})  with $y^l$  implies that
\begin{gather}\label{TP47}
\mu'_ih_{jk}+\mu'_jh_{ki}+\mu'_kh_{ij}=0.
\end{gather}
By contracting  (\ref{TP47}) with $g^{jk}$, we conclude the
following
\begin{gather*}
(n+1)\mu'_i=0.
\end{gather*}
Then on a stretch Finsler spaces, $\mu_i$ is constant along any
geodesics.
\end{proof}

\pdfbookmark[1]{References}{ref}
\LastPageEnding


\begin{thebibliography}{99}

\footnotesize\itemsep=0pt


\bibitem{AZ}
Akbar-Zadeh H.,
 Sur les espaces de Finsler \`a courbures sectionnelles constantes,
{\it Acad. Roy. Belg. Bull. Cl. Sci. (5)}
{\bf 74}  (1988), no.~10, 271--322.


\bibitem{An}
Antonelli P.L.,
Handbook of Finsler geometry, Kluwer Academic Publishers, Dordrecht, 2003.


\bibitem{Be}
 Berwald L.,
 \"{U}ber Parallel\"{u}bertragung in R\"{a}umen mit allgemeiner Massbestimmung,
 {\it Jahresbericht D.M.V.} {\bf 34} (1926), 213--220.


 \bibitem{BM}
 B\'{a}cs\'{o} S., Matsumoto M.,
 On Finsler spaces of Douglas type~-- a generalization of notion of Berwald space,
 {\it Publ. Math. Debrecen} {\bf 51} (1997), 385--406.


\bibitem{BP}
 B\'{a}cs\'{o} S., Papp I.,
 A note on a generalized  Douglas space,
 \href{http://dx.doi.org/10.1023/B:MAHU.0000038974.24588.83}{{\it Period. Math. Hungar.}} {\bf 48} (2004), 181--184.


\bibitem{ChSh}
Chen X., Shen Z.,
On Douglas metrics,
{\it Publ. Math. Debrecen} {\bf 66} (2005), 503--512.


\bibitem{Ich}
Ichijy\={o} Y.,
Finsler manifolds modeled on a Minkowski space,
{\it J. Math. Kyoto Univ.} {\bf 16} (1976), 639--652.


\bibitem{Mat}
Matsumoto M.,
On $C$-reducible Finsler spaces,
{\it Tensor (N.S.)} {\bf 24}  (1972), 29--37.


 \bibitem{MH}
 Matsumoto M., H\={o}j\={o} S.,
 A conclusive theorem for $C$-reducible Finsler spaces,
 {\it Tensor (N.S.)} {\bf 32}  (1978), 225--230.


\bibitem{NST1}
 Najaf\/i B., Shen Z., Tayebi A.,
 On a projective class   of Finsler metrics,
 {\it Publ. Math. Debrecen} {\bf 70}  (2007), 211--219.


\bibitem{NST2}
 Najaf\/i B., Shen Z., Tayebi A.,
 Finsler metrics of scalar f\/lag curvature with special non-Riemannian curvature properties,
\href{http://dx.doi.org/10.1007/s10711-007-9218-9}{{\it Geom. Dedicata}} {\bf 131}  (2008), 87--97.


\bibitem{PTP}
 Pande H.D., Tripathi P.N., Prasad B.N.,
 On a special form of the $hv$-curvature tensor of Berwald's connection $B\Gamma$ of Finsler space,
 {\it Indian J. Pure. Appl. Math.} {\bf 25}  (1994), 1275--1280.


\bibitem{ShLec}
Shen Z.,
Lectures on Finsler geometry, World Scientif\/ic Publishing Co., Singapore, 2001.


\bibitem{ShDiff}
 Shen Z.,
 Dif\/ferential geometry of spray and Finsler spaces, Kluwer Academic Publishers, Dordrecht, 2001.


\bibitem{TAE}
 Tayebi A., Azizpour E., Esraf\/ilian E.,
 On a family of connections in Finsler geometry,
 {\it Publ. Math. Debrecen} {\bf 72} (2008), 1--15.


\bibitem{TR}
 Tayebi A., Raf\/ie Rad M.,
S-curvature of isotropic Berwald metrics,
\href{http://dx.doi.org/10.1007/s11425-008-0095-y}{{\it Sci. China Ser. A}} {\bf 51}   (2008), 2198--2204.

\end{thebibliography}
\end{document}